\documentclass[11pt]{amsart}

\newtheorem{theorem}{Theorem}[section]
\newtheorem{lem}[theorem]{Lemma}
\newtheorem{prop}[theorem]{Proposition}

\theoremstyle{definition}
\newtheorem{definition}[theorem]{Definition}
\newtheorem{example}[theorem]{Example}

\theoremstyle{remark}

\numberwithin{equation}{section}

\newcommand{\ra}{\rightarrow}

\newcommand{\ddt}{\frac{\text{d}}{\text{dt}}}

\newcommand{\lb}{\langle}
\newcommand{\rb}{\rangle}

\newcommand{\mg}{\mathfrak{g}}
\newcommand{\mh}{\mathfrak{h}}
\newcommand{\mk}{\mathfrak{k}}

\renewcommand{\mp}{\mathfrak{p}}

\newcommand{\T}{{\rm T}}
\newcommand{\R}{\mathbb{R}}

\newcommand{\ii}{\mathbf{i}}
\newcommand{\jj}{\mathbf{j}}
\newcommand{\kk}{\mathbf{k}}

\begin{document}

\newcommand{\spacing}[1]{\renewcommand{\baselinestretch}{#1}\large\normalsize}
\spacing{1.14}

\title[Invariant metrics]{Invariant metrics with nonnegative curvature on $SO(4)$ and other Lie groups}

\author{Jack Huizenga, Kristopher Tapp}
\address{Department of Mathematics\\University of Chicago, Chicago, IL 60637}
\email{huizenga@uchicago.edu}
\address{Department of Mathematics\\ Williams College\\
Williamstown, MA 01267}
\email{ktapp@williams.edu}
\thanks{Supported in part by NSF grants DMS--0303326 and DMS-0353634.}
\date{\today}

\begin{abstract}
We develop techniques for classifying the nonnegatively curved
left-invariant metrics on a compact Lie group $G$.  We prove
rigidity theorems for general $G$ and a partial classification for
$G=SO(4)$.  Our approach is to reduce the general question to an
infinitesimal version; namely, to classify the directions one can
move away from a fixed bi-invariant metric such that curvature
variation formulas predict nearby metrics are nonnnegatively curved.
\end{abstract}
\maketitle

\section{Introduction}
The starting point for constructing all known examples of compact
manifolds with positive (or even quasi-positive) curvature is the
fact that bi-invariant metrics on compact Lie groups are
nonnegatively curved.  In order to generalize this fundamental
starting point, we address the question: given a compact Lie group
$G$, classify the left-invariant metrics on $G$ which have
nonnegative curvature.  New examples could potentially, via familiar
quotient constructions, lead to new examples of quasi-positively
curved spaces.  On the other hand, proofs that there are no new
examples would serve as further evidence that the known
constructions are rigid and canonical.

The first two cases, $G=SO(3)$ and $U(2)$, were completely solved
in~\cite{SMALL}.  For $G=U(2)$, all such metrics lie in the closure
of those coming from Cheeger's method, which is essentially the only
known construction of nonnegatively curved left-invariant metrics.
These classifications made use of techniques that only work in low
dimensions.  For higher-dimensional groups, more tools are necessary
to approach the problem effectively.  One important new tool is the
following, which implies in particular that the nonnegatively curved
metrics form a path-connected subset within the space of all
left-invariant metrics.

\begin{theorem}\label{TILC}
If $h$ is a left-invariant metric with nonnegative curvature on a
compact Lie group $G$, then the unique inverse-linear path from any
fixed bi-invariant metric $h(0)$ to $h(1)=h$ is through
nonnegatively curved metrics.
\end{theorem}

Here, a path of inner products on $\mg=T_eG$ (or the induced path of
left-invariant metrics) is called \emph{inverse-linear} if the the
inverses of the associated path of symmetric matrices form a
straight line.  So to classify the left-invariant metrics on $G$
with nonnegative curvature, we can first classify the directions
$h'(0)$ one can go away from a fixed bi-invariant metric $h(0)$ such
that the inverse-linear path $h(t)$ appears (up to derivative
information at $t=0$) to remain nonnegatively curved. Then, for each
candidate direction, we must check how far nonnegative curvature is
maintained along that path.

This is the approach we use for general $G$.  In the case $G=SO(4)$,
our results provide strong evidence that all left-invariant metrics
lie in the closure of those coming from Cheeger's method; that is,
there do not seem to be any new examples.  One of our stronger
results towards the classification for $SO(4)$ is the following.

\begin{theorem}\label{SO4 Main Theorem}
If $h$ is a left-invariant metric with nonnegative curvature on
$SO(4)$ and if the matrix of $h$ has an eigenvector in one of the
simple factors of $so(4)=so(3)\oplus so(3)$, then $h$ is a known
example of a metric of nonnegative curvature.
\end{theorem}

The known examples come from Cheeger's method via an action of
$T^2$ or $S^3$, as explained in Section 7.  Those from a $T^2$
action have a singular eigenvector, as in the above theorem.

 The authors are pleased to thank Burkhard Wilking, Craig
Sutton, Emily Proctor, Zachary Madden, Nela Vukmirovic, Angela
Doyle, Min Kim and the referee for numerous helpful discussions
and comments on this work.

\section{Cheeger's method}
In this section, we review Cheeger's method for altering a
nonnegatively curved metric via a group of isometries, and use it to
prove Theorem~\ref{TILC}.

Let $(M,h_0)$ be a nonnegatively curved manifold on which a compact
Lie group $G$ acts by isometries.  Let $h_R$ be a right-invariant
metric on $G$ with nonnegative curvature (often chosen to be
bi-invariant).  Notice that $G$ acts on $M\times G$  as
$g\star(p,a)=(g\star p,ag^{-1})$.  The orbit space is diffeomorphic
to $M$ via the map $[p,g]\mapsto g\star p$.  Consider the
one-parameter family of induced nonnegatively curved Riemannian
submersion metrics, $h_t$, on this orbit space:
$$(M,h_t) = \left(M\times\left(G,(1/t)h_R\right)\right)/G.$$
This family extends smoothly at $t=0$ to the original metric $h_0$
on $M$.  To describe the metric variation at a fixed $p\in M$, let
$\{v_1,...,v_k\}\subset T_p M$ denote the values at $p$ of the
Killing fields on $M$ associated to an $h_R$-orthonormal basis
$\{e_1,...,e_k\}$ of the Lie algebra $\mg$ of $G$.  Cheeger's
formula in~\cite{Cheeger} implies that the path of matrices
$A^t_{ij}=h_t(v_i,v_j)$ evolves according to
\begin{equation}\label{CH}A^t = A^0(I+tA^0)^{-1}.\end{equation}

Several authors have derived curvature-variation formulas, although
they usually assume $h_R$ is bi-invariant;
see~\cite{Muter},\cite{Strake},\cite{PW},\cite{Wil}. For this, it is
useful to consider the bijection $\Phi_t:T_p M\ra T_p M$ which
describes $h_t$ in terms of $h_0$ in the sense that for all $X,Y\in
T_p M$,
$$h_t(X,Y) = h_0(\Phi_t(X),Y).$$
This family of inner products on $T_pM$ is \emph{inverse-linear}.
This means that the path $t\mapsto \Phi_t^{-1}$ is linear, so
$\Phi_t=(I-t\Psi)^{-1}$ for some endomorphism $\Psi:\T_p M\ra T_p
M$.

Cheeger mentioned that $h_t$
has no more zero-curvature planes than $h_0$.  A precise formulation of this comment, found for example in~\cite{PW}, is
\begin{lem}\label{cheegermonotone} If the plane $\sigma=\text{span}\{X,Y\}$ has positive curvature with
respect to $h_0$, then the plane
$\Phi_t^{-1}(\sigma)=\text{span}\{\Phi_t^{-1}(X),\Phi_t^{-1}(Y)\}$
has positive curvature with respect to $h_t$.
\end{lem}
So the most natural
variational approach is to differentiate the curvature with respect
to $h_t$ of the plane $\Phi_t^{-1}(\sigma)$; this was systematically
studied in~\cite{Muter}.  In the next section, we will borrow and
generalize this idea.

\begin{proof}[Proof of Theorem~\ref{TILC}]
Let $h$ be a left-invariant metric with nonnegative curvature on the
compact Lie group $G$.  Let $h_0$ be a fixed bi-invariant metric on
$G$.  Consider the family $h_t$ of nonnegatively curved metrics on
$G$ defined by
$$(G,h_t) = ((G,h_0)\times(G,(1/t)h))/G,$$
where $G$ acts diagonally on the right of both factors.  For this
action to be isometric, $h$ must be re-considered as a
right-invariant metric on $G$, which is no problem because the left-
and right-invariant metrics determined by an inner product on $\mg$
are isometric via the inversion map.  Notice that each $h_t$ is a
left-invariant metric on $G$.

Let $\{E_1,...,E_k\}$ be an $h_0$-orthonormal basis of $\mg$ which
diagonalizes $h$.  Let $\{\lambda_1,...,\lambda_k\}$ be the
corresponding eigenvalues of $h$, so that
$\{e_i=E_i/\sqrt{\lambda_i}\}$ is an $h$-orthonormal basis of $\mg$.
In Formula~\ref{CH}, $v_i=e_i$ and $A^0=\text{diag}(1/\lambda_i)$,
so $A^t = \text{diag}(1/(\lambda_i+t))$. Thus, in the basis
$\{E_i\}$, the matrix for $\Phi_t$ is
$$\Phi_t=\text{diag}(1+(1/\lambda_i)t)^{-1}.$$
Therefore, $\Phi_t=(I-t\Psi)^{-1}$, where
$\Psi=\text{diag}(-1/\lambda_i)$.  We see that, as previously
mentioned, the path is inverse-linear.

There is no value of $t$ for which $h_t = h$.  Instead we will show
that the path $h_t$ (for $t\in[0,\infty)$) visits scalings of all of
the metrics along the unique inverse-linear path $\tilde{h}_s$
between $\tilde{h}_0=h_0$ and $\tilde{h}_1=h$.  Let $\tilde{\Phi}_s$
determine this path, so that $\tilde{h}_s(X,Y)=h_0(\tilde{\Phi}_s
X,Y)$ for all $X,Y\in\mg$.  We have that
$\tilde{\Phi}_s=(I-s\tilde{\Psi})^{-1}$, where $\tilde{\Psi}$ with
respect to the basis $\{E_i\}$ is given by
$$\tilde{\Psi}=I-\tilde{\Phi}_1^{-1} = \text{diag}(1-1/\lambda_i).$$
It is easy to see that the paths $\tilde{\Phi}_s$ (for $s\in[0,1)$)
and $\Phi_t$ (for $t\in[0,\infty)$) visit the same family of metrics
up to scaling.  More precisely, $c\cdot\tilde{\Phi}_s = \Phi_t$ when
$t=s/(1-s)$ and $c=1-s$.
\end{proof}

The method of the proof can be used to connect any two nonnegatively
curved left-invariant metrics $h_1$ and $h_2$ on $G$ through a path
of nonnegatively curved metrics.  The resultant path of inner
products on $\mg$ is inverse-linear, but this is largely irrelevant
to the question at hand because the path is not through
left-invariant metrics.


\section{Curvature Variation of Zero-Planes}\label{intro}
In this and the next section, we derive a curvature-variation
formula for an inverse-linear path of left-invariant metrics
beginning at a bi-invariant metric.

Let $G$ be a compact Lie group.  Let $h_t$ be an inverse-linear path
of left-invariant metrics on $G$ beginning at a bi-invariant metric
$h_0$.  The value of $h_t$ at $e$ is determined in terms of $h_0$ by
some self-adjoint $\Phi_t:\mg\ra\mg$ defined so that for all
$X,Y\in\mg$,
$$h(X,Y)=h_0(\Phi_t(X),Y).$$
Recall that ``inverse-linear'' means that
$$\Phi_t=(I-t\Psi)^{-1}$$ for some endomorphism $\Psi:\mg\ra\mg$.  Notice that
$\Psi=\frac{d}{dt}|_{t=0}\Phi_t$, and therefore $\Psi$ is
$h_0$-self-adjoint.  For fixed $X,Y\in\mg$, define $\kappa(t)$ to be
the unnormalized sectional curvature of
$\{\Phi_t^{-1}X,\Phi_t^{-1}Y\}$ with respect to the metric $h_t$.
The domain of $\kappa(t)$ is the open interval of $t$'s for which
$\Phi_t$ represents a nondegenerate metric; this interval depends on
the eigenvalues of $\Psi$.

Two important decisions here are inspired by properties of Cheeger's
method: (1) restricting to inverse-linear paths, and (2)
``twisting'' the plane whose curvature we are tracking.  Even though
we are considering general paths, not necessarily arising from
Cheeger's method, Theorem~\ref{TILC} and several results to follow
indicate that these decisions provide the correct approach.

If $Z_1,$ $Z_2\in \mg$, we write $\langle Z_1,Z_2\rangle =
h_0(Z_1,Z_2)$, $|Z_1|^2 = h_0(Z_1,Z_1)$, and $|Z_1|_{h_t}^2 =
h_t(Z_1,Z_1)=\langle \Phi_t Z_1,Z_1\rangle$.  We first describe $\kappa(t)$ in the important special
case where $[X,Y]=0$, so that $\kappa(0)=(1/4)|[X,Y]|^2=0$. In other words, we
first study the variation of curvature for an initially zero
curvature plane.

\begin{prop}\label{P:com} If $[X,Y]=0$, then $\kappa(0)=0$, $\kappa'(0)=0$, $\kappa''(0)=0$ and
\begin{eqnarray*}
(1/6)\kappa'''(0) & = & \lb[X,\Psi Y]+[\Psi X,Y],[\Psi X,\Psi Y]\rb
                   +\lb[\Psi X,X],\Psi[\Psi Y,Y]\rb \\
             &   &   -\lb[X,\Psi Y],\Psi[X,\Psi Y]\rb
                   -\lb[X,\Psi Y],\Psi[\Psi X,Y]\rb
                   -\lb[\Psi X,Y],\Psi [\Psi X,Y]\rb.
\end{eqnarray*}
Moreover, for all $t$ in the domain of $\kappa$,
$$\kappa(t) = t^3\cdot(1/6)\kappa'''(0)
  - t^4\cdot(3/4)|[\Psi X,\Psi Y]-\Psi([\Psi X,Y]+[X,\Psi Y])|^2_{h_t}$$
\end{prop}

We will prove this proposition in the next section as a special case
of a more general formula which does not assume that $X$ and $Y$
commute.

In the Taylor series of $\kappa(t)$ at $0$, the first non-vanishing
derivative is the third, after which the remaining tail sums to a
nonpositive term involving the norm with respect to $h_t$ of the
vector
$$D=[\Psi X,\Psi Y]-\Psi([\Psi X,Y]+[X,\Psi Y]).$$
In light of our formula for $\kappa(t)$, we can make the following
definition.
\begin{definition}\label{2parts}
We call $\Psi$ (or the variation $\Phi_t$) \emph{infinitesimally
nonnegative} if the following equivalent conditions hold:
\begin{enumerate} \item
For all $X,Y\in \mg$, there exists $\epsilon>0$ such that
$\kappa(t)\geq 0$ for $t\in [0,\epsilon)$.
\item For all commuting pairs $X,Y\in\mg$,
$\kappa'''(0)\geq 0$, and $\kappa'''(0)=0$ implies that $D=0$.
\end{enumerate}
\end{definition}

If in the first condition a single choice of $\epsilon>0$ works for
all pairs $X,Y$, then $\Phi_t$ has nonnegative curvature for
$t\in[0,\epsilon)$.  In this case, we call the variation
\emph{locally nonnegative}.  We do not know if infinitesimally
nonnegative implies locally nonnegative.  In any case, the infinitesimally nonnegative $\Psi$ are the candidate directions; the best available derivative information predicts that the paths in these directions are through nonnegatively curved metrics.

It is significant that the tail of the power series for
$\kappa(t)$ is nonpositive.  In addition to demonstrating the
equivalence of the two parts of Definition~\ref{2parts}, this
nonpositivity property immediately implies the following weak
version of Theorem \ref{TILC}: if $h_t$ is nonnegatively curved
for some $t>0$, then $\Psi$ is infinitesimally nonnegative.  This is the only version of Theorem~\ref{TILC} we will need throughout the rest of the paper.  It says that one will locate all nonnegatively curved metrics by searching only along the infinitesimally nonnegative paths.

If one omits the plane twisting and instead defines $\kappa(t)$ as
the unnormalized sectional curvature of $\{X,Y\}$, then
$\kappa(0)=0$ implies $\kappa'(0)=0$ and $\kappa''(0) = |[X,\Psi
Y]+[\Psi X,Y]|^2$.  This is true without assuming the path is
inverse-linear, so long as $\Psi=\frac{d}{dt}|_{t=0}\Phi_t$.  It is
interesting that $\kappa''(0)\geq 0$, but because of this, the
untwisted set-up provides little help in deciding which variations
remain nonnegatively curved.  We will stick with the twisted version
for the remainder of the paper.

\begin{example}\label{EE} Suppose $H\subset G$ is a Lie subgroup with Lie algebra $\mh\subset\mg$.  For $A\in\mg$, let $A^\mh$ and $A^{\mp}$ denote the projections of $A$ onto and orthogonal to $\mh$ with respect to $h_0$.  The variation $\Phi_t(A)=\frac{1}{1+t}A^{\mh} + A^{\mp}$ is inverse-linear and has nonnegative curvature for $t>0$.  In this variation, vectors tangent to $H$ are gradually shrunk.  The parametrization looks natural when re-described as a family of submersions metrics: $(G,h_t)=((G,h_0)\times(H,(1/t)h_0))/H$.  The $t=0$ derivative is $\Psi A=-A^{\mh}$.  Proposition~\ref{P:com} yields:\begin{equation}\label{shrink}(1/6)\kappa'''(0)=|[X^\mh,Y^\mh]|^2.\end{equation}
Equation~\ref{shrink} (together with Lemma~\ref{cheegermonotone} and the nonpositivity of the tail of the power series for $\kappa(t)$) re-proves Eschenburg's formula from~\cite{Eschenburg}, which says that with respect to the metric $h_t$ (for fixed $t>0)$), the plane spanned by $\Phi^{-1}_t(X)$ and $\Phi^{-1}_t(Y)$ has zero-curvature if and only if $[X,Y]=0$ and $[X^\mh,Y^\mh]=0$.

The full domain of this variation is $(-1,\infty)$.  As $t$ decreases from zero towards $-1$, vectors tangent to $H$ are enlarged.  Considering negative values of $t$ for this variation is equivalent to considering positive values of $t$ for the variation in the opposite direction, $-\Psi$.  For this opposite variation, $(1/6)\kappa'''(0)=-|[X^\mh,Y^\mh]|^2.$  So expanding $\mh$ immediately creates some negative curvature unless $[X^\mh,Y^\mh]=0$ whenever $[X,Y]=0$.  If $\mh$ is abelian, then $\kappa'''(0)=0$ for all commuting $X,Y$, which suggests that enlarging an abelian subalgebra might preserve nonnegative curvature.  Indeed, it is proven in~\cite{GZ} that enlarging an abelian subalgebra as far as $4/3$ always preserves nonnegative curvature.  In Section~\ref{S:enlarge}, we will study this variation in greater depth to determine which subalgebras can be enlarged without losing nonnegative curvature.
\end{example}

Notice that for $a>0$, $\Psi$ and $a\Psi$ generate different parameterizations of the same family of metrics.  A slightly less obvious equivalence involves adding a multiple of the identity to $\Psi$.
\begin{prop} If $\Psi$ is infinitesimally nonnegative, then so is $\tilde{\Psi}=\Psi+a\cdot I$ for any $a>0$.
\end{prop}
This proposition gives the correct equivalence modulo which one should classify the infinitesimally nonnegative endomorphisms $\Psi$.
\begin{proof} $\Psi$ and $\tilde{\Psi}$ yield the same values for $\kappa'''(0)$ and $D$ in Proposition~\ref{P:com}.  To verify this, it is convenient to use Equation~\ref{gamma0}.

An alternative proof is to observe that the inverse-linear paths $\Phi(t)=(I-t\Psi)^{-1}$ and $\tilde{\Phi}(s)=(I-s\tilde{\Psi})^{-1}$ visit the same family of metrics, modulo scalings and re-parameterizations.  More precisely, $c\cdot \Phi(t) =\tilde{\Phi}(s)$ as long as $c = 1-s\cdot a$ and $t = s/(1-s\cdot a)$.  Notice this idea was used previously in the proof of Theorem~\ref{TILC}.
\end{proof}


\section{Curvature Variation of general planes}\label{S:power
series}
In this section we state and prove a generalization of
Proposition~\ref{P:com} which does not assume $X$ and $Y$ commute.
We use this result to prove the proposition.

Certain elements of $\mg$ will appear
frequently in what follows, so to simplify the exposition we
introduce the Lie algebra elements
\begin{eqnarray*}
A &=& [\Psi X,Y]+[X,\Psi Y]\\
B &=& [\Psi X,\Psi Y]\\
C &=& [\Psi X,Y]+[\Psi Y,X]\\
D &=& \Psi^2 [X,Y] -\Psi A +B.
\end{eqnarray*}
The definition of $D$ given here coincides with the definition of
the previous section when $X$ and $Y$ commute.

\begin{theorem}\label{kappa formula theorem}
For any $t$ in the domain of $\kappa$, \begin{equation}\label{kappa
formula}\kappa(t) = \alpha + \beta t + \gamma t^2 + \delta t^3 -
\frac 34 t^4\cdot |D|_{h_t}^2,\end{equation} where
\begin{eqnarray*}
\alpha &=& \frac{1}{4} |[X,Y]|^2\\
\beta &=& -\frac{3}{4} \langle \Psi [X,Y],[X,Y]\rangle\\
\gamma &=& -\frac 34
|\Psi[X,Y]|^2+\frac 32 \langle \Psi [X,Y],A\rangle-\frac 12 \langle [X,Y],B\rangle \\
&&-\frac 14 |A|^2 +\frac 14 |C|^2 - \langle [\Psi X,X],[\Psi
Y,Y]\rangle\\
 \delta &=& -\frac 34
 \langle \Psi^3 [X,Y],[X,Y]\rangle+\frac 32 \langle \Psi^2 [X,Y],A
 \rangle - \frac 32 \langle \Psi
 [X,Y],B\rangle   \\&& - \frac 34 \langle \Psi A,A\rangle-\frac 14 \langle \Psi C,
 C\rangle+\langle \Psi[\Psi X,X], [\Psi Y,Y]\rangle+\langle
 A,B\rangle.
\end{eqnarray*}
\end{theorem}

There are two steps to the proof of this theorem.  First we prove
that Equation \ref{kappa formula} holds for all sufficiently small
$t$.  Next we show that each side of the equation is analytic.  This
allows us to invoke the well-known identity theorem: if $f,$
$g\colon I\to \R$ are analytic on an open interval $I$ and $f$ and
$g$ agree on a subinterval of $I$, then $f=g$.  We therefore
conclude that Equation \ref{kappa formula} holds for all $t$. To
accomplish the first step, we calculate the Taylor series of
$\kappa(t)$ at $t=0$. This calculation will also serve as the
foundation for our analyticity arguments.

\begin{prop}\label{Power Series}
The Taylor series of $\kappa(t)$ at $0$ is given by $$\kappa(t) =
\alpha+\beta t + \gamma t^2 +\delta t^3 -\frac 34 \sum_{n=4}^\infty
t^n \langle \Psi^{n-4} D,D\rangle,$$ with convergence for $|t|<
\|\Psi\|^{-1}$, where $\|\Psi\|=\sup_{|X|=1} |\Psi X|$ is the
operator norm of $\Psi$.
\end{prop}
\begin{proof}
In \cite{put}, P\"uttmann shows that the unnormalized sectional
curvature of vectors $Z_1,$ $Z_2\in \mg$ with respect to a
left-invariant metric $h$ whose matrix with respect to $h_0$ is
$\Phi$ is given by
\begin{eqnarray}
k_{h}(Z_1,Z_2) &=& \frac 12 \langle [\Phi Z_1,Z_2] + [Z_1,\Phi
Z_2],[Z_1,Z_2]\rangle - \frac 34 |[Z_1,Z_2]|_h^2 \nonumber \\
&&+\frac 14 \langle [Z_1,\Phi Z_2] +[Z_2,\Phi Z_1],\Phi^{-1}
([Z_1,\Phi Z_2] + [Z_2,\Phi Z_1])\rangle \label{Puttmanns formula}\\
&&-\langle [Z_1,\Phi Z_1], \Phi^{-1} [Z_2,\Phi Z_2]\rangle.\nonumber
\end{eqnarray} It follows that
\begin{eqnarray*}
\kappa(t) &=& k_{h_t}(\Phi_t^{-1} X,\Phi_t^{-1} Y)\\
&=& \frac 12 \langle [X,\Phi_t^{-1} Y] +[\Phi_t^{-1}
X,Y],[\Phi_t^{-1}X,\Phi_t^{-1}Y]\rangle\\
&&-\frac 34 \langle \Phi_t [\Phi_t^{-1} X,\Phi_t^{-1}Y],[\Phi_t^{-1}
X,\Phi_t^{-1} Y]\rangle\\
&&+\frac 14 \langle [\Phi_t^{-1} X,Y] + [\Phi_t^{-1}
Y,X],\Phi_t^{-1} ([\Phi_t^{-1} X,Y]+[\Phi_t^{-1} Y,X])\rangle\\
&& - \langle [\Phi_t^{-1} X,X],\Phi_t^{-1} [ \Phi_t^{-1}
Y,Y]\rangle\\
&=& I_1 - I_2 + I_3 - I_4.
\end{eqnarray*}
Using the expression $\Phi_t^{-1} = I-t\Psi$, we can easily simplify
$I_1$, $I_3$, and $I_4$.  We find
\begin{eqnarray*}
I_1 &=& |[X,Y]|^2 - \frac {3t}2 \langle [X,Y],A\rangle + t^2(\langle
[X,Y],B\rangle +\frac 12 |A|^2 \rangle - \frac {t^3}2 \langle
A,B\rangle\\
I_3 &=& \frac {t^2}4 |C|^2 - \frac {t^3}4\langle C,\Psi C\rangle\\
I_4 &=& t^2 \langle [\Psi X,X],[\Psi Y,Y]\rangle - t^3 \langle [\Psi
X,X],\Psi[\Psi Y,Y]\rangle.
\end{eqnarray*}
To calculate $I_2$, notice that if $|t|< \|\Psi\|^{-1}$, then
$$\Phi_t = \sum_{n=0}^\infty t^n \Psi^n,$$ with convergence
in the space of endomorphisms of $\mg$ with the operator norm. From
this formula we calculate
\begin{eqnarray*}
\frac 43 I_2 &=&  \langle \Phi_t([X,Y]-t A +t^2 B),[X,Y]-t A +t^2
B\rangle\\
&=&  \sum_{n=0}^\infty t^n \langle \Psi^n [X,Y] - t \Psi^n A
+ t^2 \Psi^n B,[X,Y]-t A +t^2 B\rangle\\
&=&  \sum_{n=0}^\infty t^n \left(\langle \Psi^n[X,Y],[X,Y]\rangle-2t
\langle \Psi^n [X,Y],A\rangle \right.\\&&\left.+t^2(\langle \Psi^n
A,A\rangle+2\langle \Psi^n [X,Y],B\rangle)-2t^3 \langle \Psi^n
A,B\rangle+t^4 \langle \Psi^n B,B\rangle \right)\\
&=& |[X,Y]|^2 +t(\langle \Psi [X,Y],[X,Y]\rangle  - 2 \langle
[X,Y],A\rangle)\\&&+t^2(\langle \Psi^2 [X,Y],[X,Y]\rangle-2\langle
\Psi [X,Y],A\rangle+|A|^2+2\langle [X,Y],B\rangle)\\
&&+t^3(\langle \Psi^3[X,Y],[X,Y]\rangle - 2 \langle
\Psi^2[X,Y],A\rangle +\langle \Psi A,A\rangle \\&&\quad+2\langle
\Psi [X,Y],B\rangle-2\langle A,B\rangle)
\\&&+\sum_{n=4}^\infty t^n \langle \Psi^{n-4} D,D\rangle.\end{eqnarray*}
Combining the different terms proves the result.
\end{proof}

Notice the power series of $\kappa(t)$ would have been much messier
if we were considering the unnormalized sectional curvature of $X$
and $Y$ with respect to $h_t$ instead of the unnormalized sectional
curvature of $\Phi_t^{-1} X$ and $\Phi_t^{-1} Y$. The value of
twisting is even apparent at a purely computational level.

When $|t|<\|\Psi\|^{-1}$, we observe $$ -\frac 34 \sum_{n=4}^\infty
t^n\langle \Psi^{n-4} D,D\rangle  = -\frac 34 t^4 \langle \Phi_t
D,D\rangle = -\frac 34 t^4 \cdot |D|_{h_t}^2.$$ This proves Equation
\ref{kappa formula} holds for small $t$. Therefore to complete the
proof of Theorem \ref{kappa formula theorem}, all we must do is
prove $\kappa(t)$ and $|D|_{h_t}^2$ are analytic.

\begin{lem}\label{analyticity}
The function $\kappa(t)$ is analytic on its domain of definition.
\end{lem}
\begin{proof}
Assume that $t_0$ is such that $\Phi_{t_0}$ corresponds to a metric
on $G$. We show $\kappa$ is locally a power series at $t_0$.
Recalling P\"utmann's Formula \ref{Puttmanns formula}, it is clear
we must only prove that
$$|[\Phi_t^{-1} X,\Phi_t^{-1} Y]|^2_{h_t}$$ can be expressed as a power
series near $t_0$.  Since $\Psi$ is $h_0$-self-adjoint, it can be
diagonalized; say $\Psi=\text{diag}(a_1,\ldots,a_d)$.  We then have
\begin{eqnarray}\label{Phi t equation} \Phi_t&=&\text{diag}\left(\frac
1{1-a_1t},\ldots,\frac 1{1-a_dt}\right)\nonumber\\
&=& \text{diag}\left(\frac{1}{1-a_it_0} \sum_{n=0}^\infty
\left(\frac{a_i}{1-a_i t_0}\right)^n
(t-t_0)^n\right)\\
&=& \Phi_{t_0} \sum_{n=0}^\infty \Phi_{t_0}^n \Psi^n
(t-t_0)^n\nonumber,\end{eqnarray} with convergence whenever
$|t-t_0|$ is sufficiently small.  We can use this expression for
$\Phi_t$ together with the identity $\Phi_t^{-1} = I-t_0 \Psi
-(t-t_0)\Psi$ to expand $|[\Phi_t^{-1} X,\Phi_t^{-1} Y]|_{h_t}^2$ as
a power series as in the proof of Proposition \ref{Power Series}.
\end{proof}

Analyticity of $|D|_{h_t}^2$ also follows from Equation \ref{Phi t
equation}, completing the proof of Theorem \ref{kappa formula
theorem}.

\begin{proof}[Proof of Proposition \ref{P:com}]
Assume $X$ and $Y$ commute.  It is easy to see $\alpha=\beta=0$, and
that $\delta$ equals $6$ times the stated formula for
$\kappa'''(0)$. All that remains to be shown is $\gamma=0$.  But the
bi-invariance of $h_0$ and the Jacobi identity give the identity
\begin{eqnarray}\label{gamma0} \langle[\Psi
X,Y],[X,\Psi Y]\rangle &=& - \langle\Psi X,[[X,\Psi Y],Y]\rangle =
\langle\Psi X,[[\Psi Y,Y],X]+[[Y,X],\Psi Y]\rangle
\\\nonumber &=& \langle\Psi X,[[\Psi Y,Y],X]\rangle = - \langle[\Psi X,X],[\Psi Y,Y]\rangle,
\end{eqnarray}
from which $\gamma=0$ follows easily.
\end{proof}


\section{A general rigidity result}

The next lemma is our primary tool for deriving rigidity statements
about infinitesimally nonnegative variations; it plays an important
role in Section 7, where we give a partial classification of the
infinitesimally nonnegative endomorphisms of $so(4)$.

\begin{lem}\label{k}
Assume that $\Psi$ is infinitesimally nonnegative. Let $\mp_0$ be
the eigenspace of $\Psi$ corresponding to the smallest eigenvalue.
If $X\in\mp_0$, $Y\in\mg$ and $[X,Y]=0$, then $[X,\Psi Y]\in\mp_0$.
\end{lem}
\begin{proof}
Proposition \ref{P:com} applied to $X$ and $Y$ gives:
$$(1/6)\kappa'''(0) = a_0|[X,\Psi Y]|^2 - \lb[X,\Psi Y],\Psi [X,\Psi Y]\rb,$$
where $a_0$ is the smallest eigenvalue.  This is negative unless
$[X,\Psi Y]\in\mp_0$.
\end{proof}

The next proposition is a global version of this lemma.  The
argument used in its proof serves as the prototype for how we
transform rigidity statements about infinitesimally nonnegative
endomorphisms into rigidity statements about nonnegatively curved
metrics.

\begin{prop}\label{global rigidity}
Assume that $\Phi$ is the matrix of a nonnegatively curved metric,
$h$. Let $\mp_0$ be the eigenspace of $\Phi$ corresponding to the
smallest eigenvalue.  If $X\in \mp_0$, $Y\in \mg$ and $[X,Y]=0$,
then $[X,\Phi^{-1} Y]\in \mp_0$.
\end{prop}
\begin{proof}
Let $\Psi=I-\Phi^{-1}$, so that $\Phi_t=(I-t\Psi)^{-1}$ is the
unique inverse-linear path from $h_0$ to $h_1=h$.  Theorem~\ref{TILC}
says $\Psi$ must be infinitesimally nonnegative.  Notice that $\Psi$ and $\Phi$ have
the same smallest eigenspace $\mp_0$.  Proposition~\ref{k} gives
that $$[X,\Psi Y] = [X,(I-\Phi^{-1})Y] = -[X,\Phi^{-1}Y]\in\mp_0.$$
\end{proof}

We note that this result can also be derived directly from
P\"uttmann's Formula \ref{Puttmanns formula}.


\section{Enlarging subalgebras}\label{S:enlarge}

Here we continue the discussion on enlarging subalgebras begun in
Example \ref{EE}.  Let $H\subset G$ be a Lie subgroup of the Lie
group $G$ with Lie algebra $\mh\subset \mg$.  For $Z\in \mg$, denote
by $Z^\mh$ and $Z^\mp$ the projections of $Z$ onto $\mh$ and its
$h_0$-orthogonal complement $\mp$.  Let $\Psi(Z)=Z^\mh$, so $\Phi_t
= (I-t\Psi)^{-1}$ is the inverse-linear variation which gradually
expands vectors in $\mh$ as $t$ increases from $0$. If $\mh$ is
abelian, it is easy to use the formulas for the coefficients of the
power series of $\kappa(t)$ in tandem with the analyticity of
$\kappa$ to prove
\begin{equation}\label{Abelian blowup curv eq}\kappa(t) = \frac 14 |[X,Y]|^2 - \frac 34 |[X,Y]^\mh|^2
\cdot \frac t{1-t} \qquad (-\infty < t < 1).\end{equation} From this
formula we can show that enlarging $\mh$ by a factor of up to $4/3$
always preserves nonnegative curvature, a result which first
appeared in \cite{GZ}.  In fact, the particularly nice form of
$\kappa(t)$ allows us to prove a stronger statement.

\begin{theorem}\label{Abelian Blowup}
Scaling the abelian subalgebra $\mh\subset \mg$ preserves
nonnegative curvature if and only if no vector in $[\mg,\mg]$ has
the square of its norm expanded by more than $4/3$.
\end{theorem}
\begin{proof}
By Equation \ref{Abelian blowup curv eq}, the metric $h_t$ is
nonnegatively curved if and only if
\begin{equation}\label{blowup inequality} |Z^\mh|^2 \cdot \frac t {1-t} \leq \frac 13|Z|^2
 \end{equation} holds for all $Z\in [\mg,\mg]$. As
 $$|Z|_{h_t}^2 = \langle \Phi_t Z,Z\rangle = \langle Z + \frac
 t{1-t} Z^\mh,Z\rangle = |Z|^2 + |Z^\mh|^2 \cdot \frac t{1-t},$$ we
 find Inequality \ref{blowup inequality} is equivalent to
 requiring that $|Z|_{h_t}^2 \leq (4/3)\cdot |Z|^2$ holds for all
 $Z\in [\mg,\mg]$.
\end{proof}

If $[\mg,\mg]\cap \mh\neq \{0\}$, this theorem says that $\mh$ can be scaled up by a
factor up to $4/3$. At the other extreme, if $[\mg,\mg]\perp \mh$ then
we find that $\mh$ can be expanded up by an arbitrary amount.  This
was already known, since if $\mh$ is orthogonal to $[\mg,\mg]$ then
$\mh$ is contained in the center of $\mg$.  This rescaling then stays
within the family of bi-invariant metrics on $\mg$.

When $\mh$ is not abelian, things are not quite so simple.  In this
case the power series simplifies to
$$\kappa(t) = \frac 14 |[X,Y]|^2 -\frac 34 |[X,Y]^\mh|^2 t +\frac 34
|B|^2 t^2 - \frac 14 |B|^2 t^3 -\frac 34 |[X^\mp,Y^\mp]^\mh|^2 \cdot
\frac {t^2}{1-t}.$$ We can use this formula to classify exactly
which subalgebras of $\mg$ can be enlarged a small amount while
maintaining nonnegative curvature.

\begin{theorem}\label{Arbitrary Blowup}
Expanding the subalgebra $\mh\subset \mg$ by a small amount preserves
nonnegative curvature if and only if there exists a constant $c$
such that $|[X^\mh,Y^\mh]|\leq c\cdot |[X,Y]|$ holds for all $X,$
$Y\in \mg$.
\end{theorem}

We omit the lengthy but easy proof for the reason that we do not
know if there are any interesting examples of subalgebras for which
the latter condition holds.  It clearly holds when $\mh$ is either
abelian or an ideal of $\mg$ (or the sum of an ideal and an
orthogonal abelian subalgebra), but it is already known that such
subalgebras can be enlarged while maintaining nonnegative curvature.


\section{Known metrics on $SO(4)$ with nonnegative curvature}
Each known example of a left-invariant metric $h$ with nonnegative
curvature on $G=SO(4)$ comes from Cheeger's construction.  In this
section, we catalog each known example in terms of the eigenvalue
and eigenvector structure of the map $\Phi$ representing it with
respect to a fixed bi-invariant metric $h_0$, meaning that
$h(A,B)=h_0(\Phi A,B)$.

\subsection{Product Metrics} The Lie algebra $\mg=so(4)$
is a product $\mg=\mg_1\oplus\mg_2$, with each factor isomorphic to
$so(3)$.  The two factors are $h_0$-orthogonal.  If they are
$h$-orthogonal, then $h$ is a product metric on $SO(4)$'s double
cover $S^3\times S^3$.  The classification of product metrics with
nonnegative curvature reduces to the classification of
left-invariant metrics with nonnegative curvature on $SO(3)$, solved
in~\cite{SMALL}.  Observe that for any product metric, $\mg$
decomposes into three 2-dimensional $\Phi$-invariant abelian
subalgebras, obtained by pairing eigenvectors from the two factors.

As for infinitesimal examples, if $\Psi$ is a product, meaning
$\Psi(\mg_1)\subset\mg_1$ or equivalently $\Psi(\mg_2)\subset\mg_2$,
then the inverse-linear path $\Phi_t=(I-t\Psi)^{-1}$ it generates is
though product metrics, which have nonnegative curvature for small
$t$.

\subsection{Torus Actions} Let $\{A_1,A_2,A_3\}$ and $\{B_1,B_2,B_3\}$ be
$h_0$-orthonormal bases of $\mg_1$ and $\mg_2$, respectively.  After
scaling $\mg_1$ and $\mg_2$ by factors $c$ and $d$, respectively,
then enlarging the abelian subalgebra $\tau=\text{span}\{A_3,B_1\}$
by $4/3$, then further altering the metric on $\tau$ via the
remaining $T^2$-action on $G$, one obtains a nonnegatively curved
metric $h$ with matrix $\Phi$ of the form
\begin{equation}\label{Torus Form}\left(\begin{array}{cccccc} c & 0 & 0 & 0 & 0 & 0\\
0 & c & 0 & 0 & 0&0\\
0 & 0 &a_1 & a_3 & 0 & 0\\
0 & 0 &a_3 & a_2 & 0 & 0\\
0 & 0 & 0 & 0 & d &0\\
0 &0 &0 & 0 & 0 & d
\end{array}\right)\end{equation}
with respect to the basis $\{A_1,A_2,A_3,B_1,B_2,B_3\}$.  In the
final alteration, any right-invariant (and hence bi-invariant and
flat) metric on $T^2$ can be used.  The only restriction on $\Phi$,
coming from the fact that this final alteration only shrinks
vectors, is that the norm on $\tau$ determined by the matrix
$\left(\begin{array}{cc} a_1 & a_3 \\ a_3 & a_2\end{array}\right)$
is strictly bounded above by the norm determined by
$\left(\begin{array}{cc} \frac 43 \cdot c & 0 \\ 0 & \frac 43 \cdot
d\end{array}\right).$  Limit points of such metric are also
nonnegatively curved.  That is, we must consider the closure of the
known examples, which transforms the strict inequality above into a
non-strict one.

Observe that $\mg$ decomposes into three 2-dimensional
$\Phi$-invariant abelian subalgebras: one equals $\tau$, and the
other two are obtained by pairing vectors in $\mg_1$ with vectors in
$\mg_2$.

Notice that any endomorphism $\Psi$ with the matrix form of
Equation~\ref{Torus Form} will generate an inverse-linear variation
$\Phi_t=(I-t\Psi)^{-1}$.  These metrics will be nonnegatively curved
for some interval $t\in[0,\epsilon)$.  The parameters
$\{c,d,a_1,a_2,a_3\}$ defining $\Psi$ are unrestricted, although
they do determine $\epsilon$.

\subsection{$S^3$-actions} Let $\tilde{h}$ denote the bi-invariant metric on
$S^3\times S^3$ obtained from $h_0$ by rescaling $\mg_1$ and $\mg_2$
by factors $a$ and $b$ respectively. Let $g_R$ denote a
right-invariant metric with nonnegative curvature on $S^3$ with
eigenvalues $\{\lambda_1,\lambda_2,\lambda_3\}$ and eigenvectors
$\{e_1,e_2,e_3\}$.  Define a metric $h$ by
$$(S^3\times S^3,h)=((S^3\times S^3,\tilde{h})\times(S^3,g_R))/S^3,$$
where $S^3$ acts diagonally.  Consider the basis
$$\mg=\mg_1\oplus\mg_2=\text{span}\{A_1,A_2,A_3\}\oplus\text{span}\{B_1,B_2,B_3\},$$
where $A_i=(e_i,0)$ and $B_i=(0,e_i)$.  Let
$V_i=\text{span}\{A_i,B_i\}$, which for each $i$ is a 2-dimensional
abelian subalgebra of $\mg$. Notice that the three $V_i$'s are
mutually orthogonal with respect to $h_0$, $\tilde{h}$, and $h$.  It
therefore suffices to describe $h$ in terms of $h_0$ separately on
each $V_i$.

For this, the matrix representing $\tilde{h}$ in terms of $h_0$ on
$V_i$ in the basis $\{A_i,B_i\}$ is $M_i=\left(\begin{matrix} a & 0
\\ 0 & b\end{matrix}\right)$.  The matrix representing $h$ in terms
of $\tilde{h}$ in the basis $\{A_i+B_i,b A_i - a B_i\}$ is
$N_i=\left(\begin{matrix} t_i & 0 \\ 0 & 1\end{matrix}\right)$,
where $t_i=\frac{\lambda_i}{1+\lambda_i}$. Thus, letting $T$ be the
change of basis matrix, $T=\left(\begin{matrix} 1 & b \\ 1 &
-a\end{matrix}\right)$, the matrix we seek which represents $h$ in
terms of $h_0$ on $V_i$ in the basis $\{A_i,B_i\}$ is
\begin{equation}\label{willow}\Phi_i= M_i(TN_iT^{-1})=\frac{1}{a+b}
\left(\begin{matrix} a(b+at_i) & ab(t_i-1) \\ ab(t_i-1) & b(a+bt_i)\end{matrix}\right).\end{equation}

In summary, $\mg$ decomposes into the three $\Phi$-invariant
2-dimensional abelian subalgebras, $\{V_1,V_2,V_3\}$.  However, with
only the five parameters $\{a,b,t_1,t_2,t_3\}$ under our control,
and with restrictions on the $t$'s, we do not attain the full
9-parameter family of metrics for which the subalgebras
$\{V_1,V_2,V_3\}$ are $\Phi$-invariant.

Infinitesimal examples have the form $\Psi:=I-\Phi^{-1}$ with $\Phi$
in the form of Equation~\ref{willow}.  A calculation shows that all
such matrices have the form
$\Psi=\text{diag}(\Psi_1,\Psi_2,\Psi_3)$, where
\begin{equation}\label{gg}\Psi_i=\left(\begin{matrix}\alpha & 0 \\ 0 & \beta\end{matrix}\right) -
\frac{1}{2\lambda_i}\left(\begin{matrix} 1 & 1 \\ 1 &
1\end{matrix}\right).\end{equation} The parameters $\alpha,\beta$
are free, but the parameters $\{\lambda_1,\lambda_2,\lambda_3\}$ are
restricted to be eigenvalues of a nonnegatively curved metric on
$SO(3)$.


\section{Infinitesimal rigidity for $SO(4)$}
In this section, we assume that $G=SO(4)$ and $\Psi:\mg\ra\mg$ is
infinitesimally nonnegative, and we prove rigidity results for
$\Psi$.  In the next section, we translate these infinitesimal
rigidity results into global theorems.

Recall that $\mg = so(4)=\mg_1\oplus\mg_2$ is a product, and
$X\in\mg$ is called \emph{regular} if it has non-zero projections
onto both $\mg_1$ and $\mg_2$; otherwise, it is called
\emph{singular}. We give $G$ the most natural bi-invariant metric $h_0$,
so that any orthonormal bases of the factors $\mg_1$ and $\mg_2$
behave like the quaternions $\{\ii,\jj,\kk\}$ with respect to their
Lie bracket structure.  We will show in Section 10 that there is no
essential loss of information in restricting ourselves to only
working with this bi-invariant metric.

The previous section classified the known possibilities of $\Psi$
into three types, coming from: (1) products, (2) torus actions and
(3) $S^3$-actions.  In the first two cases, $\Psi$ has a non-zero
singular eigenvector, while in the third case, it does not.

\begin{theorem}\label{Th1} If $\Psi$ has a non-zero singular eigenvector, the either $\Psi$ is
a product or $\Psi$ has the form of Equation~\ref{Torus Form}.  In
either case, $h_t$ is a family of known examples with nonnegative
curvature for sufficiently small $t$.\end{theorem}

If $\Psi$ has no non-zero singular eigenvectors, we conjecture that
$\Psi$ is a known example coming from an $S^3$-action.  A first step
in this direction is to locate three $\Psi$-invariant abelian
subalgebras.  The following theorem falls just short of this goal:
\begin{theorem}\label{Th2}
There are orthonormal bases $\{A_1,A_2,A_3\}$ and $\{B_1,B_2,B_3\}$
of the two factors of $\mg=\mg_1\oplus\mg_2$ such that with respect
to the basis $\{A_1,B_1,A_2,B_2,A_3,B_3\}$, $\Psi$ has the form
$$\Psi=\left(\begin{matrix}
a_1 & a_3 & 0 & 0 & 0 & 0 \\
a_3 & a_2 & 0 & 0 & 0 & 0 \\
0 & 0 & b_1 & b_3 & \lambda & 0\\
0 & 0 & b_3 & b_2 & 0 & \mu\\
0 & 0 & \lambda & 0 & c_1 & c_3 \\
0 & 0 & 0 & \mu & c_3 & c_2
\end{matrix}\right).$$
\end{theorem}
We conjecture that $\lambda=\mu=0$, which means that $\mg$
decomposes into three orthogonal $\Psi$-invariant abelian
subalgebras, as it should.  Even granting this conjecture, there
remains the work of reducing the above 9-parameter family to the
5-parameter family of known examples from Equation~\ref{gg}.  This
appears to be a computationally difficult problem.

The remainder of this chapter is devoted to proving
Theorems~\ref{Th1} and~\ref{Th2}.  We begin with a weak version of
Theorem~\ref{Th1}.  Recall that $\mp_0$ denotes the eigenspace
corresponding to the smallest eigenvalue, $a_0$, of $\Psi$.
\begin{lem}\label{nongen}
If $\mp_0$ contains a non-zero singular vector, then either $\Psi$
is a product or $\Psi$ has the form of Equation~\ref{Torus Form}.
\end{lem}
\begin{proof}
Without loss of generality, assume there exists a non-zero vector
$X_1\in\mg_1\cap\mp_0$.  Assume that $\Psi$ is not a product, so
there exists $\hat{Y}\in\mg_2$ such that $\Psi\hat{Y}$ has a nonzero
projection, $X_2$, onto $\mg_1$.  Notice that $X_1$ and $X_2$ are
orthogonal because
$$\lb X_1,X_2\rb = \lb X_1,\Psi \hat{Y}\rb = \lb \Psi X_1,\hat{Y}\rb = a_0\lb X_1,\hat{Y}\rb = 0.$$
Let $X_3=[X_1,\Psi\hat{Y}]\in\mg_1$, which by Lemma~\ref{k} lies in
$\mp_0$, so $\text{span}\{X_1,X_3\}\subset\mp_0$.  Let $Y_2$ be the
projection of $\Psi X_2$ onto $\mg_2$, which is a non-zero vector by
the self-adjoint property of $\Psi$. Complete $\{Y_2\}$ to an
orthogonal basis $\{Y_1,Y_2,Y_3\}$ of $\mg_2$, ordered so that their
bracket structure is like $\{\ii,\jj,\kk\}$. Notice that $\Psi
(\text{span}\{Y_1,Y_3\})\subset\mg_2$ (again by the self-adjoint
property of $\Psi $).  In summary, after scaling all the vectors to
unit-length, we have an orthonormal basis:
$$\mg = \mg_1\oplus\mg_2 = \text{span}\{X_1,X_2,X_3\}\oplus\text{span}\{Y_1,Y_2,Y_3\}$$
with $\text{span}\{X_1,X_3\}\subset\mp_0$, and $\Psi  X_2 =
cY_2+\lambda X_2$ (for some $c,\lambda\in\R$ with $c\neq 0$), and
$\Psi (\text{span}\{Y_1,Y_3\})\subset\mg_2$.

Applying Proposition~\ref{P:com} to the vectors $X_2$ and $Y_1$
gives
\begin{eqnarray*}
\kappa'''(0) & = & 6\lb[\Psi  X_2,Y_1],[\Psi X_2,\Psi Y_1]\rb
             - 6\lb [\Psi  X_2,Y_1],\Psi [\Psi  X_2,Y_1]\rb \\
        & = & 6\lb [ cY_2,Y_1],[cY_2,\Psi  Y_1]\rb - 6\lb[cY_2,Y_1],\Psi [cY_2,Y_1]\rb\\
        & = & -6c^2\lb Y_3,[Y_2,\Psi  Y_1]\rb - 6c^2\lb Y_3,\Psi  Y_3\rb \geq 0.
\end{eqnarray*}
Notice that
\begin{eqnarray*}
\lb Y_3,[Y_2,\Psi  Y_1]\rb & = & \lb Y_3,[Y_2,\text{projection of }\Psi Y_1\text{ onto }Y_1]\rb \\
                             & = & \lb Y_3,[Y_2,\lb\Psi Y_1,Y_1\rb Y_1]\rb \\
                             & = & -\lb \Psi Y_1,Y_1\rb,
\end{eqnarray*}
from which we conclude
$$\lb Y_1,\Psi Y_1\rb\geq\lb Y_3,\Psi Y_3\rb.$$
Similarly, applying Proposition \ref{P:com} to the vectors $X_2$ and
$Y_3$ yields the reverse inequality, so:
$$\lb Y_1,\Psi Y_1\rb = \lb Y_3,\Psi Y_3\rb.$$

Replacing $Y_1$ and $Y_3$ with any other orthonormal basis of
$\text{span}\{Y_1,Y_3\}$ yields the same conclusion.  In other
words, for any angle $\theta$, if we set $a=\cos(\theta)$ and
$b=\sin(\theta)$ then
$$\lb aY_1+bY_3,\Psi (aY_1+bY_3)\rb = \lb bY_1-aY_3,\Psi (bY_1-aY_3)\rb.$$
This implies that $\lb Y_1,\Psi Y_3\rb = \lb \Psi Y_1,Y_3\rb = 0$.
The linear map from $\text{span}\{Y_1,Y_3\}$ to $\R$ sending
$Y\mapsto\lb \Psi Y,Y_2\rb$ has a non-zero vector in its kernel.
Assume without loss of generality that $Y_1$ is in its kernel.
Notice that $Y_1$ is an eigenvector of $\Psi $.

In the ordered basis $\{X_1,X_2,X_3,Y_1,Y_2,Y_3\}$, we thus far have
$$\Psi =\left(\begin{matrix}
a_0 & 0 & 0 & 0 & 0 & 0 \\
0 & \lambda & 0 & 0 & c & 0 \\
0 & 0 & a_0 & 0 & 0 & 0\\
0 & 0 & 0 & \beta & 0 & 0\\
0 & c & 0 & 0 & \gamma & s \\
0 & 0 & 0 & 0 & s & \beta
\end{matrix}\right)$$
Applying our $\kappa'''(0)$ formula to $X=X_2$ and $Y=aY_2+bY_3$
gives
$$\kappa'''(0) = 6bc^2(as+b\beta) - 6 b^2c^2\beta = 6bc^2as.$$
Since $\kappa'''(0)\geq 0$ for all choices of $\{a,b\}$, and $c\neq
0$, we learn that $s=0$.  After re-ordering the basis, $\Psi$ has
the form of Equation~\ref{Torus Form}.
\end{proof}
\begin{theorem}\label{franklin}
The eigenspace $\mp_0$ contains a non-zero vector which belongs to a
$\Psi$-invariant $2$-dimensional abelian subalgebra of $\mg$.
\end{theorem}
\begin{proof}
If $\mp_0$ contains a non-zero singular vector, the conclusion
follows easily from Lemma~\ref{nongen}, so we assume that this is
not the case.  When $A=(A_1,A_2)\in\mg=\mg_1\oplus\mg_2$ is regular,
let
$\overline{A}=\left(\frac{|A_2|}{|A_1|}A_1,-\frac{|A_1|}{|A_2|}A_2\right)$,
which commutes with $A$, is orthogonal to $A$, and has the same norm
as $A$.

The proof is indirect.  We assume for each $A\in\mp_0$ that
$\text{span}\{A,\overline{A}\}$ is not $\Psi$-invariant, and we
derive a contradiction.

Let $A\in\mp_0$ be unit-length.  Since $\Psi $ is self-adjoint,
$\Psi \overline{A}$ is orthogonal to $A$.  Notice that
$\overline{A}$ is not an eigenvector of $\Psi $; if it were, then
$\text{span}\{A,\overline{A}\}$ would be an invariant abelian
subalgebra.  Therefore, $[A,\Psi \overline{A}]$ is non-zero.  Let
$B$ be the unit-length vector in the direction of $[A,\Psi
\overline{A}]$.  By Lemma~\ref{k}, $B\in\mp_0$.  Notice that $B$ is
orthogonal to $A$ and $\overline{A}$.

So far we know that $\dim(\mp_0)\geq 2$.  Clearly $\dim(\mp_0)\leq
3$ because it contains no non-zero singular vectors, and hence
intersects $\mg_1$ and $\mg_2$ trivially.  We wish to prove
$\dim(\mp_0)=2$.  Suppose to the contrary that $\dim(\mp_0)=3$.
Consider the map from $\mp_0$ to $\mp_0$ defined as
$$Z\mapsto [Z,\Psi\overline{Z}].$$
By the above arguments, this map sends each unit-length $Z\in\mp_0$
to a non-zero vector in $\mp_0$ orthogonal to $Z$.  This map
therefore induces a smooth non-vanishing vector field on the unit
$2$-sphere in $\mp_0$, which is a contradiction.  Thus,
$\dim(\mp_0)=2$.  Notice $A$ and $B$ play symmetric roles in that
$[B,\Psi\overline{B}]$ is parallel to $A$ (because it lies in
$\mp_0$ and is perpendicular to $B$), and $A$ is orthogonal to $B$
and $\overline{B}$.

Choose unit-length vectors $C_1\in\mg_1$ and $C_2\in\mg_2$ such that
$\{A,\overline{A},B,\overline{B},C_1,C_2\}$ is an orthonormal basis
of $\mg$.  For $i=1,2$, the $\mg_i$-components of $\{A,B,C_i\}$ form
an orthogonal basis of $\mg_i$.  The $C_i$'s can be chosen so that
these orthogonal bases are oriented, so after normalizing, they act
like $\{\ii,\jj,\kk\}$ with respect to their Lie bracket structure.
For purposes of calculating Lie brackets in this basis, we lose no
generality in assuming that for some $a,b\in(0,1)$,
\begin{gather}\label{iii}
A=(a\ii,\sqrt{1-a^2}\ii), \,\,\,\, B=(b\jj,\sqrt{1-b^2}\jj),\,\,\,\, C_1=(\kk,0) \\
\overline{A}=(\sqrt{1-a^2}\ii,-a\ii), \,\,\,\, \overline{B}=(\sqrt{1-b^2}\jj,-b\jj),\,\,\,\, C_2=(0,\kk).\notag
\end{gather}
Notice that $\langle \Psi\overline{A},\overline{B}\rangle=\langle
\Psi\overline{B},\overline{A}\rangle=0$, because if
$\Psi\overline{A}$ had a nonzero $\overline{B}$-component, then
$[A,\Psi\overline{A}]$ would have nonzero $C_1$ and
$C_2$-components.

In the basis $\{A,\overline{A},B,\overline{B},C_1,C_2\}$, $\Psi $
has the form
\begin{equation}\label{early}\Psi =\left(\begin{matrix}
a_0 & 0 & 0 & 0 & 0 & 0 \\
0 & p & 0 & 0 & \alpha_1 & \alpha_2 \\
0 & 0 & a_0 & 0 & 0 & 0\\
0 & 0 & 0 & q & \beta_1 & \beta_2\\
0 & \alpha_1 & 0 & \beta_1 & f_1 & f_2 \\
0 & \alpha_2 & 0 & \beta_2 & f_2 & f_3
\end{matrix}\right).\end{equation}
There are a few obvious restrictions among the variables determining
$\Psi $.  For example, since $[A,\Psi \overline{A}]$ is parallel to
$B$, and $[B,\Psi \overline{B}]$ is parallel to $A$, we learn
\begin{equation}\label{simon}\frac{\alpha_1}{\alpha_2}=\frac{\beta_2}{\beta_1}
=\frac{b\sqrt{1-a^2}}{a\sqrt{1-b^2}},\end{equation} and we obtain
\begin{equation}\label{early2}\Psi =\left(\begin{matrix}
a_0 & 0 & 0 & 0 & 0 & 0 \\
0 & p & 0 & 0 & \alpha & \alpha \cdot s \\
0 & 0 & a_0 & 0 & 0 & 0\\
0 & 0 & 0 & q & \beta \cdot s & \beta\\
0 & \alpha & 0 & \beta \cdot s & f_1 & f_2 \\
0 & \alpha \cdot s & 0 & \beta & f_2 & f_3
\end{matrix}\right),\end{equation}
where $s = \frac{ a \sqrt{1-b^{2}} }{ b \sqrt{1-a^{2}} } > 0$ and
$\alpha, \beta \neq 0$.

Using Lemma~\ref{k}, we can now prove that $s=1$ and consequently $a
= b$. Indeed,  for every $Z\in\text{span}\{A,B\}$, we have $[Z,\Psi
\overline{Z}]\in\text{span}\{A,B\}$.  In particular, let $Z_t=(\cos
t) A+(\sin t)B$, so
$$\overline{Z}_t = \left(  f(t) \left(a\cos(t) \ii + b\sin(t) \jj\right),  - (1/f(t)) \left( \sqrt{1-a^{2}} \cos(t) \ii + \sqrt{ 1-b^{2}}\sin(t)\jj \right)\right),$$
where
$$f(t)=  \sqrt{\frac{ (1-a^{2})\cos^{2}(t) + (1-b^{2})\sin^{2}(t) }{ a^{2}\cos^{2}(t) + b^{2}\sin^{2}(t)} }.$$
We will use that the following vector lies in $\text{span}\{A,B\}$:
\begin{eqnarray*}
Q &=& \ddt\Big|_{t=0}[Z_t,\Psi\overline{Z_t}] = [B,\Psi\overline{A}] + \left[A, \Psi \left(f'(0)a\ii + f(0) b \jj, -g'(0) \sqrt{1-a^{2}}  \ii -g(0) \sqrt{ 1-b^{2}} \jj\right)\right] \\
    &=& [B,\Psi\overline{A}] + \left[A, \Psi \left(f(0) b \jj, -g(0) \sqrt{ 1-b^{2}} \jj\right)\right] \\
    &=& [B,\Psi\overline{A}] + \left[A, \Psi \left( \frac{b \sqrt{1-a^{2}}}{a} \jj , -\frac{a \sqrt{1-b^{2}}}{\sqrt{ 1-a^{2} }} \jj \right)\right] \\
    &=& [B,\Psi\overline{A}] + \left[A, \Psi \left( \sqrt{1-b^{2}} \cdot  s^{-1} \, \jj , -b \cdot  s \, \jj\right)\right].
\end{eqnarray*}
In particular, $Q$ is perpendicular to $\overline{A}$, so
\begin{eqnarray*}
0 &=& \lb Q,\overline{A}\rb = \lb [B,\Psi\overline{A}],\overline{A}\rb +
\left\langle\left[A, \Psi \left( \sqrt{1-b^{2}} \cdot  s^{-1} \, \jj , -b \cdot  s \, \jj\right)\right],\overline{A}\right\rangle \\
   &=&  \lb [B,\Psi\overline{A}],\overline{A}\rb = - \langle \Psi \overline{A}, [B, \overline{A} ] \rangle \\
   &=& - \langle p \overline{A} + (\alpha\kk , \alpha s\kk) , [(b\jj, \sqrt{1-b^{2}} \jj), (\sqrt{1-a^{2}} \ii, -a\ii)] \rangle \\
   &=& - \langle p \overline{A} + (\alpha\kk , \alpha s\kk) , ( - b \sqrt{1-a^{2}} \kk, a \sqrt{1-b^{2}} \kk) \rangle \\
   &=& \alpha b \sqrt{1-a^{2}} - s\alpha a \sqrt{1-b^{2}},
\end{eqnarray*}
which implies $s = \frac{b \sqrt{1-a^{2}}}{a \sqrt{1-b^{2}} } =
s^{-1}$. It follows that $s=1$ and, consequently, $a= b$. Now the
fact that the orthogonal projection of $Q$ onto
$\text{span}\{C_1,C_2\}$ is zero is equivalent to
\begin{equation}\label{gab}p(-b\sqrt{1-a^2}\kk,a\sqrt{1-b^2}\kk)+q(a\sqrt{1-b^2}\kk,-b\sqrt{1-a^2}\kk)=0.\end{equation}
Since $a=b$, this implies that $q=p$. So we obtain
\begin{equation}\label{psi}\Psi =\left(\begin{matrix}
a_0 & 0 & 0 & 0 & 0 & 0 \\
0 & p & 0 & 0 & \alpha & \alpha \\
0 & 0 & a_0 & 0 & 0 & 0\\
0 & 0 & 0 & p & \beta & \beta\\
0 & \alpha & 0 & \beta & f_1 & f_2 \\
0 & \alpha & 0 & \beta & f_2 & f_3
\end{matrix}\right).\end{equation}
Since $a=b$, it is easy to see that
$[A,\overline{B}]+[B,\overline{A}]=0$. This implies
$V_1=\beta\overline{A}-\alpha\overline{B}$ commutes with $V_2=\beta
A-\alpha B.$  Since $V_2\in\mp_0$, and $V_1$ is an eigenvector of
$\Psi$ (with eigenvalue $p$), we learn that $\text{span}\{V_1,V_2\}$
is a $\Psi$-invariant 2-dimensional abelian subalgebra of $\mg$
containing a non-zero vector in $\mp_0$.  This is a contradiction.
\end{proof}

\begin{proof}[Proof of Theorem~\ref{Th2}]
By the previous theorem, there exists a $\Psi$-invariant abelian
subalgebra of $\mg$, spanned by some $A_1\in\mg_1$ and some
$B_1\in\mg_2$.  Let $V_1$ denote the orthogonal compliment of $A_1$
in $\mg_1$, and let $V_2$ denote the orthogonal compliment of $B_1$
in $\mg_2$.

Let $\pi_1:\mg\ra\mg_1$  and $\pi_2:\mg\ra\mg_2$ denote the
projections.  Define $T_1:V_1\ra V_2$ as
$T_1=\pi_2\circ\Psi|_{V_1}$, and define $T_2:V_2\ra V_1$ as
$T_2=\pi_1\circ\Psi|_{V_2}$.  Notice that for all $A\in V_1$ and
$B\in V_2$,
$$\lb T_1A,B\rb = \lb\Psi A,B\rb = \lb A,\Psi B\rb = \lb A,T_2 B\rb.$$

Let $S^1$ denote the circle of unit-length vectors in $V_1$.  Let
$R:S^1\ra S^1$ denote a $90^\circ$ rotation.  Define $F:S^1\ra\R$ by
$F(A)=\lb T_1(A),T_1(R(A))\rb$.  For all $A\in S^1$,
$$F(R(A))=\lb T_1(R(A)),T_1(-A)\rb = -F(A).$$
This implies that there exists $A_2\in S^1$ such that $F(A_2)=0$.
Let $A_3=R(A_2)$. First suppose $T_1$ (and hence also $T_2$) is
nonsingular. Define $B_2=T_1(A_2)/|T_1(A_2)|$ and
$B_3=T_1(A_3)/|T_1(A_3)|$.  The fact that $F(A_2)=0$ immediately
implies $B_2$ and $B_3$ are orthogonal, and that $T_2(B_2)\parallel
A_2$ and $T_2(B_3)\parallel A_3$.  Thus, the basis
$\{A_1,A_2,A_3,B_1,B_2,B_3\}$ satisfies the conclusion of the
theorem.

If $T_1$ (and hence also $T_2$) is singular, then arbitrary
orthonormal bases $\{A_2,A_3\}$ of $V_1$ and $\{B_2,B_3\}$ of $V_2$
work, so long as $A_2\in\text{ker}(T_1)$ and
$B_2\in\text{ker}(T_2)$.
\end{proof}
Our final proof in this section is due to Nela Vukmirovic and
Zachary Madden:
\begin{proof}[Proof of Theorem~\ref{Th1}]
Choose bases $\{A_1,A_2,A_3\}$ of $\mg_1$ and $\{B_1,B_2,B_3\}$ of
$\mg_2$ so that $\Psi$ has the matrix form of Theorem~\ref{Th2}.
With respect to the ordering $\{A_3,A_2,A_1,B_1,B_2,B_3\}$, $\Psi$
then has the form $$\label{psi}\Psi =\left(\begin{matrix}
c_1 & \lambda & 0 & 0 & 0 & c_3\\
\lambda & b_1 & 0 & 0 & b_3 & 0 \\
0 & 0 & a_1 & a_3 & 0 & 0\\
0 & 0 & a_3 & a_2 & 0 & 0\\
0 & b_3 & 0 & 0 & b_2 & \mu \\
c_3 & 0 & 0 & 0 & \mu & c_2
\end{matrix}\right).$$ If $a_3=0$, then the result follows from
Lemma~\ref{nongen}, so we can assume $a_3\neq 0$. To complete the
proof, we show that $c_1=b_1$, $b_2=c_2$, and
$\lambda=\mu=b_3=c_3=0$, which puts $\Psi$ into Form \ref{Torus
Form}.  The hypothesis that $\Psi$ has a non-zero singular
eigenvector implies $b_3=0$ or $c_3$=0. Without loss of generality,
assume $b_3=0$. Henceforth, the value $\kappa'''(0)$ with respect to
the commuting pair $X=\alpha_1 A_1+\alpha_2 A_2+\alpha_3 A_3$ and
$Y=\beta_1 B_1+\beta_2 B_2+\beta_3 B_3$ will be denoted by
$[\alpha_1,\alpha_2,\alpha_3,\beta_1,\beta_2,\beta_3]$.  These
$6$-tuples are easily expanded using Maple or Mathematica.

First, $[0,\pm 1,1,1,0,0]=c_3^2(a_2-b_2)\pm 4a_3^2\lambda\geq 0$.
However, as $[0,0,1,0,1,0] + [0,0,1,0,0,1] = c_3^2(b_2-a_2) \geq 0$,
we deduce $\lambda =0$ and consequently $c_3^2(b_2-a_2)=0$.
Similarly, $[1,0,0,0,\pm 1,1] = c_3^2(a_1-b_1)\pm 4a_3^2\mu \geq 0$.
But $[1,0,0,0,1,0] + [1,0,0,0,0,1] = c_3^2(b_1-a_1) \geq 0,$ so it
follows that $\mu =0$ and $c_3^2(b_1-a_1)=0$.

Furthermore, the inequalities $[0,1,0,1,0,0]\geq 0$ and
$[0,0,1,1,0,0]\geq 0$ give respectively the plus and minus versions
of the inequality $ \pm a_3^2(b_1-c_1) \geq 0$. Analogously, from
examining $[1,0,0,0,1,0]$ and $[1,0,0,0,0,1]$ we conclude $\pm
a_3^2(b_2-c_2) \geq 0$. Since $a_3$ is non-zero we get that $b_1 =
c_1$ and $b_2=c_2$.

All that remains to be shown is that $c_3=0$.  If $c_3\neq 0$, then
$a_1=b_1$ and $a_2=b_2$. By considering $[1,1,1,1,1,1]$,
$[1,1,1,-1,1,1]$, $[1,1,1,1,-1,1]$, and $[1,1,1,1,1,-1],$ we deduce
 $\pm a_3^2c_3 \geq 0$, which implies that $c_3=0$.  Thus,
$\Psi$ has the form of Equation~\ref{Torus Form}.
\end{proof}


\section{Global rigidity for $SO(4)$}\label{SO4 Section}

The previous section partially classified the infinitesimally
nonnegative endomorphisms for $G=SO(4)$.  We now translate these
infinitesimal results into a partial classification of the
nonnegatively curved left-invariant metrics on $SO(4)$.

Assume $G=SO(4)$.  Let $\Phi$ be the matrix for a nonnegatively
curved left-invariant metric $h$ on $G$.  The variation
$\Phi_t=(I-t\Psi)^{-1}$ satisfies $\Phi_1=\Phi$ as long as we choose
$\Psi=I-\Phi^{-1}$.  By Theorem~\ref{TILC}, this variation is
through nonnegatively curved metrics, so $\Psi$ is infinitesimally
nonnegative.  We will apply restrictions on $\Psi$ from the previous
section in order to prove rigidity theorems about $\Phi$.

First, we prove a global analog of Theorem \ref{Th1}.  This theorem
implies Theorem \ref{SO4 Main Theorem} from the introduction.
\begin{theorem}\label{Singular eigenvector classification}
If $\Phi$ has a singular eigenvector, then either $h$ is a product
metric or $h$ comes from a torus action.  In either case, $h$ is a
known example of a metric of nonnegative curvature.
\end{theorem}
\begin{proof}
Since $\Phi$ has a singular eigenvector, so does $\Psi$.  According
to Theorem~\ref{Th1}, either $\Psi$ is a product or $\Psi$ can be
written in Form~\ref{Torus Form}.  If $\Psi$ is a product then
$\Phi$ is a product, which means $h$ is a product metric.  If
instead $\Psi$ has Form~\ref{Torus Form}, then so does $\Phi$.

Assume $\Phi$ has Form~\ref{Torus Form}; we must prove that $\Phi$
satisfies the $4/3$-restriction shared by all known examples.
Permuting some basis vectors if necessary, we may assume that $A_1$,
$A_2$, $A_3$ and $B_1$, $B_2$, $B_3$ behave like the quaternions
$\ii$, $\jj$, $\kk$ with respect to their Lie bracket structure.
Denote by $\tilde h$ the metric on $\tau$ corresponding to the
matrix $$ \left(\begin{array}{cc} \frac 43 \cdot c & 0 \\ 0 & \frac
43 \cdot d\end{array}\right).$$ We must prove that
\begin{equation*}|\alpha A_3+\beta B_1|_h^2
\leq |\alpha A_3 +\beta B_1|_{\tilde h}^2\end{equation*} holds for
all $\alpha,\beta \in \R$.

Consider the unnormalized sectional curvature of the vectors $\alpha
A_1+\beta B_2$ and $A_2+B_3$ with respect to $h$. We have
\begin{eqnarray*} [\Phi (\alpha A_1+\beta B_2),A_2+B_3] &=& \alpha c
A_3+\beta d B_1\\{} [\alpha A_1+\beta B_2,\Phi(A_2+B_3)] &=& \alpha
c A_3 + \beta d B_1\\{} [\alpha A_1 + \beta B_2,A_2+B_3] &=& \alpha
A_3 +\beta B_1,
\end{eqnarray*}
and therefore by P\"uttmann's Formula \ref{Puttmanns formula}
\begin{eqnarray*}
k_h(\alpha A_1+\beta B_2,A_2+B_3) &=& \langle \alpha c A_3 + \beta d
B_1,\alpha A_3+\beta B_1\rangle-\frac 34 |\alpha A_3 + \beta B_1|_h^2\\
&=& \frac 34 (|\alpha A_3+\beta B_1|_{\tilde h}^2 - |\alpha A_3
+\beta B_1|_{h}^2).
\end{eqnarray*}
Since $h$ is nonnegatively curved, this proves the required
inequality.
\end{proof}

Similarly, we obtain a global version of Theorem~\ref{Th2}.

\begin{theorem}
There are orthonormal bases $\{A_1,A_2,A_3\}$ and $\{B_1,B_2,B_3\}$
of the two factors of $\mg=\mg_1\oplus \mg_2$ such that with respect to
the basis $\{A_1,B_1,A_2,B_2,A_3,B_3\}$, $\Phi$ has the form
$$
\Phi= \left(\begin{array}{cccccc} a_1 & a_3 & 0 & 0 & 0 & 0\\ a_3 &
a_2 & 0 & 0 & 0 & 0 \\ 0 & 0 & b_1 & b_3 & \lambda & 0\\ 0 & 0 & b_3
& b_2 & 0 & \mu\\ 0 & 0 & \lambda & 0 & c_1 & c_3 \\ 0 & 0 & 0 & \mu
& c_3 & c_2 \end{array}\right).$$ In particular, $\mg$ has a
$2$-dimensional $\Phi$-invariant abelian subalgebra.
\end{theorem}
\begin{proof}
By Theorem~\ref{franklin}, $\mg$ has a 2-dimensional
$\Psi$-invariant abelian subalgebra.  This subalgebra is also
$\Phi$-invariant.  The result follows by mimicking the proof of
Theorem~\ref{Th2}.
\end{proof}


\section{Changing the initial bi-invariant metric}\label{Bi-invariant
metric shifting}

Let $h_0$ be a fixed bi-invariant metric, and consider a second
bi-invariant metric $h_1$.  If $h$ is a nonnegatively curved
left-invariant metric, then according to Theorem 1.1 the unique
inverse-linear paths from $h_0$ to $h$ and from $h_1$ to $h$ are
through nonnegatively curved metrics.  We can view this as saying
that the inverse-linear path from $h_0$ to $h$ is through
nonnegatively curved metrics if and only if the inverse-linear path
from $h_1$ to $h$ is.

In light of this result, it is natural to ask whether the
inverse-linear path from $h_0$ to $h$ is infinitesimally nonnegative
if and only if the inverse-linear path from $h_1$ to $h$ is. The
main result of this section is an affirmative answer, which shows
that the concept of ``infinitesimally nonnegative'' is independent
of the starting bi-invariant metric.  This means that when
classifying the infinitesimally nonnegative endomorphisms of $\mg$
with respect to a bi-invariant metric, the choice of bi-invariant
metric is essentially irrelevant.

\begin{theorem}\label{Inf nonneg bi-invariant metric invariance}
The inverse-linear path from $h_0$ to $h$ is infinitesimally
nonnegative if and only if the inverse-linear path from $h_1$ to $h$
is.
\end{theorem}

For the proof of this theorem, let $M$ be the matrix of $h_1$ with
respect to $h_0$, let $\Phi$ be the matrix of $h$ with respect to
$h_0$, let $\Theta$ be the matrix of $h$ with respect to $h_1$, and
put $\Psi = I-\Phi^{-1}$, $\Upsilon = I - \Theta^{-1}$.  Theorem
\ref{Inf nonneg bi-invariant metric invariance} is a consequence of
the following result.

\begin{prop}
For any commuting vectors $X$ and $Y$ in $\mg$,
$$D^{\Upsilon}_{X,Y} = D^{\Psi}_{MX,MY} \qquad \textrm{and} \qquad
\delta^{\Upsilon,h_1}_{X,Y} = \delta^{\Psi,h_0}_{MX,MY},$$ where,
for instance, $\delta_{MX,MY}^{\Upsilon,h_1}$ denotes the
coefficient $\delta$ in the power series of the function $\kappa(t)$
defined with respect to the endomorphism $\Psi$, the bi-invariant
metric $h_0$, and the commuting pair of vectors $MX,$ $MY$.  Hence
$\Psi$ is infinitesimally nonnegative if and only if $\Upsilon$ is.
\end{prop}
\begin{proof}
Write
$$\mg=\mg_1\oplus \cdots \oplus \mg_r \oplus Z(\mg),$$ where the $\mg_i$
are simple subalgebras and
$Z(\mg)$ is the center of $\mg$.  The simple subalgebras have unique
bi-invariant metrics up to a scalar multiple, any choice of inner
product on $Z(\mg)$ is bi-invariant, and all bi-invariant metrics on
$\mg$ arise as product metrics from this decomposition. We can
diagonalize $M$ with respect to a basis respecting the above
decomposition, and $M$ will have a single eigenvalue corresponding
to each simple factor $\mg_i$ and arbitrary eigenvalues on basis
vectors in $Z(\mg)$.  This allows us to factor $M=M_1\cdots M_s$,
where each $M_i$ scales an ideal of $\mg$ and leaves its orthogonal
complement fixed.  By induction, it suffices to prove the above
formulas for $M=M_1$, where $M$ acts on $\mg$ by $Z\mapsto \lambda
Z^\mh + Z^\mk$ for some $\lambda>0$ and $\mh$, $\mk$ are ideals of
$\mg$ with $\mg=\mh\oplus \mk$.  This special case follows from a
long straightforward calculation using the definitions of $D$ and
$\delta$.
\end{proof}

We conjecture that the formulas of this proposition are a special
case of a formula relating $\kappa_{X,Y}^{\Upsilon,h_1}(t)$ to
$\kappa_{MX,MY}^{\Psi,h_0}(t).$ For instance, in the special case
where $M=\lambda I$ is a scalar multiple of the identity, the
formula
$$\left(\frac{\lambda}{1-(1-\lambda)t}\right)^3 \cdot
\kappa^{\Upsilon,h_1}_{X,Y}(t) = \kappa^{\Psi,h_0}_{MX,MY}\left(
\frac{\lambda t}{1-(1-\lambda)t}\right) \qquad (0\leq t \leq 1)$$
holds, even when $X$ and $Y$ do not commute, and can be demonstrated
using the techniques of Section 4.


\bibliographystyle{amsplain}

\end{document}